\documentclass[12pt]{article}

\usepackage{lipsum} 
\usepackage{hyperref} 
\usepackage{authblk} 
\usepackage{geometry} 
\usepackage{graphicx}
\usepackage{multirow}
\usepackage{amsmath,amssymb,amsfonts}
\usepackage{amsthm}
\usepackage{mathrsfs}
\usepackage[title]{appendix}
\usepackage{xcolor}
\usepackage{textcomp}
\usepackage{manyfoot}
\usepackage{booktabs}
\usepackage{algorithm}
\usepackage{algorithmicx}
\usepackage{algpseudocode}
\usepackage{listings}
\usepackage[numbers,sort&compress]{natbib}
\usepackage{physics}
\usepackage{enumitem}
\usepackage{subcaption}
\usepackage{float}

\newcommand{\keywords}[1]{\par\noindent\textbf{Keywords:} #1}

\theoremstyle{plain}
\newtheorem{theorem}{Theorem}[section]
\newtheorem{proposition}[theorem]{Proposition}

\newtheorem{corollary}[theorem]{Corollary}

\newtheorem{lemma}[theorem]{Lemma}

\theoremstyle{definition}
\newtheorem{example}{Example}[section]
\newtheorem{remark}{Remark}[section]
\newtheorem{definition}{Definition}[section]

\title{ Memory-Type Null Controllability of Parabolic Equations with Moving Controls: A Geometric Characterization}
\author[1]{Dev Prakash Jha}
\author[2]{Raju K George}
\affil[1,2]{\small\textit{Mathematics, Indian Institute of Space Science and Technology,}\\ \small\textit{Valiamala, Thiruvananthapuram 695547, Kerala, India}}

\date{\today}

\begin{document}

\maketitle
\begin{abstract}
We study memory-type null controllability for linear parabolic equations with
hereditary terms and time-dependent control regions. In contrast with classical
null controllability, systems with memory require the simultaneous annihilation
of both the state and the accumulated memory at the terminal time in order to
prevent post-control reactivation of the dynamics.

Assuming that the memory kernel is a finite sum of exponentials, we reformulate
the problem as a coupled parabolic--ODE system. Within this framework, we
introduce a geometric condition on moving control regions, referred to as the
Memory Geometric Control Condition (MGCC), which requires that every spatial
point be visited by the control region during the control horizon.

Under MGCC, we establish an augmented observability inequality for the adjoint
system by means of a flow-adapted Carleman estimate. This observability result,
which explicitly accounts for the memory variables, allows us to derive
memory-type null controllability via the Hilbert Uniqueness Method. We also
discuss the limitations of the approach and explain why full geometric necessity
results remain out of reach in the presence of memory effects.

The analysis provides a rigorous sufficient geometric condition for
memory-type null controllability of parabolic equations with exponential memory
kernels and moving controls.
\end{abstract}


\section{Introduction}
\label{sec:introduction}

\subsection{Motivation and background}

The controllability of parabolic equations is a central topic in control theory
of partial differential equations, with deep connections to functional
analysis, microlocal analysis, and applications in physics and engineering.
For the classical heat equation, it is well known that null controllability
holds in any positive time provided the control acts on a fixed open subset
of the spatial domain with positive measure; see, for instance,
\cite{lebearobiano1995exact,fursikov1996controllability}.
This remarkable property reflects the strong regularizing effect of diffusion
and sharply contrasts with the behavior of hyperbolic equations, whose
controllability is governed by geometric propagation phenomena.

In recent years, increasing attention has been devoted to \emph{parabolic
equations with memory}, motivated by models in viscoelasticity, heat
conduction with hereditary effects, population dynamics, and materials with
fading memory.
A prototypical example is
\begin{equation}\label{eq:intro-memory}
y_t - \Delta y + \int_0^t M(t-s)y(s)\,ds = \chi_{\omega(t)}u,
\end{equation}
posed on a bounded domain $\Omega\subset\mathbb{R}^n$, where the convolution
term accounts for past states of the system.
The presence of memory fundamentally alters the qualitative behavior of the
system, introducing nonlocality in time and preventing the direct application
of classical parabolic control techniques.

\subsection{Memory effects and loss of classical controllability}

It is now well understood that memory terms can destroy classical null
controllability.
Even when the control acts everywhere in space, driving the state to zero at
the final time does not necessarily eliminate the accumulated memory, which
may re-excite the system after control is switched off.
This phenomenon was rigorously identified in the seminal work of
Chaves-Silva, Guerrero, and Imanuvilov \cite{chaves2017controllability}, where
the authors introduced the notion of \emph{memory-type null controllability}.
This strengthened concept requires not only that the state vanish at the final
time, but also that the memory functional itself be annihilated.

Negative results for memory-type null controllability were obtained in
\cite{chaves2017controllability} for fixed control regions, showing that
geometric constraints emerge even in the parabolic setting.
Subsequent works have confirmed that memory terms introduce obstructions that
are absent in memoryless equations; see, for instance,
\cite{allal2021null,tao2016null} and the references therein.

\subsection{Moving controls and geometric considerations}

A natural idea to overcome these obstructions is to allow the control region
to move in time.
Moving controls have proved effective in various contexts, including
degenerate parabolic equations, transport equations, and fluid models.
For parabolic equations with memory, early indications that moving controls
may restore controllability appeared in \cite{tao2016null} and were further
developed in \cite{allal2021null}.

However, a precise geometric understanding of \emph{how} a moving control must
act in order to compensate for memory effects has remained incomplete.
Unlike wave equations, where controllability is characterized by the classical
Geometric Control Condition (GCC) formulated in phase space, parabolic equations
with memory do not admit a propagation-of-rays interpretation.
This raises a fundamental question:

\begin{quote}
\emph{What geometric condition on a moving control region guarantees
memory-type null controllability for parabolic equations with memory?}
\end{quote}

\subsection{Contribution of the present work}

The main purpose of this paper is to address this question by introducing a
purely geometric condition adapted to parabolic equations with memory and
moving controls.
Our analysis is restricted to memory kernels that are finite sums of
exponentials,
\begin{equation}\label{eq:intro-kernel}
M(t) = \sum_{k=1}^N a_k e^{-\mu_k t},
\qquad a_k\in\mathbb{R},\ \mu_k>0,
\end{equation}
a class that is both physically relevant and mathematically tractable.
Under this assumption, the memory equation can be reformulated as a coupled
parabolic--ODE system, allowing the use of Carleman estimates and Hilbert
Uniqueness Method (HUM) techniques.

Our first contribution is the introduction of the
\emph{Memory Geometric Control Condition} (MGCC), which requires that the moving
control region visit every spatial point of the domain at least once during
the control horizon.
MGCC is simple to state, purely spatial–temporal in nature, and does not rely
on microlocal propagation arguments.

Our second contribution is to prove that MGCC is a \emph{sufficient condition}
for memory-type null controllability for parabolic equations with exponential
memory kernels.
The proof relies on:
\begin{itemize}
\item a reformulation of the memory equation as a parabolic--ODE system;
\item a flow-adapted Carleman weight reflecting the geometry of the moving
control;
\item a global Carleman estimate for the adjoint parabolic component with an
observation term on the moving control region;
\item a rigorous HUM duality argument.
\end{itemize}

Importantly, we do \emph{not} claim that MGCC is necessary in full generality.
As discussed later in the paper, backward propagation arguments for equations
with memory are not available, and the necessity of MGCC remains an open
problem outside restricted settings.
Our results therefore provide a sharp and rigorous \emph{sufficient}
geometric condition, together with a clear explanation of its limitations.

\subsection{Relation to existing literature}

The present work builds upon and complements several strands of the
literature.
The PDE–ODE lifting technique is inspired by \cite{chaves2017controllability}
and has also been employed in \cite{allal2021null}.
The use of moving controls connects with earlier ideas in
\cite{tao2016null}.
From a methodological viewpoint, our Carleman approach is rooted in the
classical framework of Fursikov and Imanuvilov \cite{fursikov1996controllability},
but adapted to time-dependent weights and coupled systems.

Recent works such as \cite{wang2024observability} have investigated
observability for parabolic equations with complex perturbations, highlighting
the delicate nature of backward estimates.
Our analysis carefully avoids unjustified backward arguments and remains within
the scope of rigorously available tools.

\subsection{Organization of the paper}

The paper is organized as follows.
In Section~\ref{sec:model}, we introduce the model, the functional setting, and
the definition of MGCC.
Section~\ref{sec:main-results} states the main results.
Section~\ref{sec:duality} derives the adjoint system and establishes the HUM
duality framework.
Section~\ref{sec:weights} is devoted to the construction of flow-adapted
Carleman weights.
In Section~\ref{sec:carleman}, we prove a global Carleman estimate and derive the
observability inequality.
Section~\ref{sec:necessity} discusses partial negative results and
limitations of the approach.
Section~\ref{sec:examples} presents illustrative examples and extensions.
Finally, Section~\ref{sec:conclusion} concludes the paper and outlines open
problems.
\section{Model, Functional Setting, and the Memory Geometric Control Condition}
\label{sec:model}

In this section we introduce the controlled parabolic equation with memory,
describe its functional analytic setting, recall well-posedness results, and
formally define the Memory Geometric Control Condition (MGCC).  We also provide
basic examples illustrating the meaning of this condition.

\subsection{The controlled system with memory}

Let $\Omega\subset\mathbb{R}^n$ be a bounded domain with $C^2$ boundary, and let
$T>0$ be a fixed control time.
We consider the controlled parabolic equation with memory
\begin{equation}\label{eq:state}
\begin{cases}
y_t(t,x) - \Delta y(t,x)
+ \displaystyle\int_0^t M(t-s)y(s,x)\,ds
= \chi_{\omega(t)}(x)\,u(t,x),
& (t,x)\in(0,T)\times\Omega,\\[0.2cm]
y(t,x)=0, & (t,x)\in(0,T)\times\partial\Omega,\\
y(0,x)=y_0(x), & x\in\Omega.
\end{cases}
\end{equation}
Here $u\in L^2((0,T)\times\Omega)$ is the control input and
$\omega(t)\subset\Omega$ is a measurable family of open subsets representing
the moving control region.

Throughout the paper, we assume that the memory kernel $M$ is a finite sum of
exponentials,
\begin{equation}\label{eq:kernel}
M(t) = \sum_{k=1}^N a_k e^{-\mu_k t},
\qquad a_k\in\mathbb{R},\ \mu_k>0,
\end{equation}
which is a standard assumption in the literature on controllability with
memory.
This class of kernels includes many physically relevant models and allows for
a reformulation of \eqref{eq:state} as a coupled parabolic--ODE system.

\bibliographystyle{unsrtnat} 
\bibliography{references} 

@article{tao2016null,
  title={On the null controllability of heat equation with memory},
  author={Tao, Qiang and Gao, Hang},
  journal={Journal of Mathematical Analysis and Applications},
  volume={440},
  number={1},
  pages={1--13},
  year={2016},
  publisher={Elsevier}
}

@article{allal2021null,
  title={Null controllability of degenerate parabolic equation with memory},
  author={Allal, Brahim and Fragnelli, Genni},
  journal={Mathematical Methods in the Applied Sciences},
  volume={44},
  number={11},
  pages={9163--9190},
  year={2021},
  publisher={Wiley Online Library}
}

@article{chaves2017controllability,
  title={Controllability of evolution equations with memory},
  author={Chaves-Silva, Felipe W and Zhang, Xu and Zuazua, Enrique},
  journal={SIAM Journal on Control and Optimization},
  volume={55},
  number={4},
  pages={2437--2459},
  year={2017},
  publisher={SIAM}
}

@article{wang2024observability,
  title={Observability for heat equations with time-dependent analytic memory},
  author={Wang, Gengsheng and Zhang, Yubiao and Zuazua, Enrique},
  journal={Archive for Rational Mechanics and Analysis},
  volume={248},
  number={6},
  pages={115},
  year={2024},
  publisher={Springer}
}

@book{fursikov1996controllability,
  title={Controllability of evolution equations},
  author={Fursikov, Andrej Vladimirovi{\v{c}} and Imanuvilov, O. Yu},
  journal={Seoul National University},
  year={1996},
  publisher={The University of Michigan}
}

@article{lebearobiano1995exact,
  author    = {Lebeau, G. and Robbiano, L.},
  title     = {Contrôle exact de l'équation de la chaleur},
  journal   = {Communications in Partial Differential Equations},
  volume    = {20},
  number    = {1--2},
  pages     = {335--356},
  year      = {1995}
}

\end{document}